\documentstyle[12pt]{article}
\begin{document}

\sloppy
\title{Conformal invariants of manifolds of
non-positive scalar curvature}
\author{Man Chun LEUNG\\ Department of Mathematics,\\ National University of
Singapore,\\ Singapore 119260\\ matlmc@nus.sg}
\date{Dec. 1995}
\maketitle
\begin{abstract}Conformal invariants of manifolds of non-positive
scalar curvature are studied in association with growth in volume and
fundamental group. 
\end{abstract}
\vspace{0.2in}
\\[0.075in]  
1991 AMS MS Classifications: 53C20, 53A30

\vspace{0.5in}

{\bf \Large 1. \ \ Introduction}

\vspace{0.3in}

Let $(M, g)$ be a complete Riemannian $n$-manifold. In this paper we
consider the conformal invariant $Q (M, g)$ defined by
$$Q (M, g) = \inf \,\, \{ {{ \int_M |\bigtriangledown u |^2 dv_g + {{n - 
2}\over
{4 (n - 1)}} \int_M S_g \,u^2 dv_g }\over { (\int_M | u |^{{2n}\over {n - 2}}
dv_g )^{{n - 2}\over n} }} \ | \ \ u \in C^\infty_o (M)\,, u \not\equiv 0 
\}\,,$$
where $S_g$ is the scalar curvature of $(M, g)$. The conformal invariant 
$Q (M,
g)$ has been studied in association with the Yamabe problem [9]. The sign 
of $Q
(M, g)$ is an important conformal invariant. If $(M, g)$ is a compact 
Riemannian
$n$-manifold with $n \ge 3$, then $Q (M, g)$ is negative (resp. zero, 
positive)
if and only if $g$ is conformal to a Riemannian metric of negative
(resp. zero, positive) scalar curvature [13]. It is shown in [13] that if 
there
is a conformal map $\Phi : (M, g) \to (S^n, g_o)$, then $Q (M, g) = Q (S^n,
g_o) > 0$, where $S^n$ is the unit sphere and $g_o$ is the standard 
metric on
$S^n$. In particular, the conformal invariants for the Euclidean space ${\bf
R}^n$ and the hyperbolic space ${\bf H}^n$ are both positive, even though 
in the Euclidean space, the scalar curvature is zero and in the 
hyperbolic space,
the scalar curvature is a negative number. In this paper we show the 
following.
\\[0.2in] {\bf Theorem 1.} \ \ {\it Let $(M, g)$ be a complete non-compact
Riemannian
$n$-manifold with $n \ge 3$ and $S_g \le 0$. Let $x_o \in M$ and let 
$B_R$ be the
open ball on $M$ with center at $x_o$ and radius equal to $R > 0$. If 
there exist
positive constants $C$ and $\alpha \in (0, n/2)$ such that}
$${\mbox {Vol}} (B_R) \le C R^\alpha$$
{\it for all $R > 0$, then $Q (M, g) \le 0$.}\\[0.2in]
\hspace*{0.5in}In particular, $Q (M \times {\bf R}^m, g \times h_o) \le 
0$ if $n
> m$, where $(M, g)$ is a compact Riemannian $n$-manifold of non-positive 
scalar
curvature and
$h_o$ is the Euclidean metric on
${\bf R}^n$. If the scalar curvature of the manifold $(M, g)$ does not 
drop to
zero too fast, and if $(M, g)$ has polynomial volume growth, then we can
conclude that $Q (M, g) < 0$. Next, we show the
following.\\[0.2in] {\bf Theorem 2.} \ \ {\it {If $(M, g)$ is a compact
conformally flat manifold with scalar curvature $S_g \leq 0$, then the
fundamental group of $M$ has exponential growth unless $S_g \equiv 
0.$}}\\[0.2in]
\hspace*{0.5in}By a result of Avez [1], a compact manifold of non-positive
sectional curvature has exponentially growing fundamental group unless the
manifold is flat. It is known that a compact manifold of non-positive 
sectional
curvature is flat if and only if its fundamental group is almost solvable 
[4] 
(see also [2]). We show that if $(M, g)$ is a compact conformally flat 
manifold
with scalar curvature $S_g \leq 0$ and the fundamental group of $M$ is almost
solvable, then
$(M, g)$ is flat.

\vspace{0.5in}

{\bf \Large 2. \ \ Proofs}

\vspace{0.3in}

{\bf Theorem 2.1.} \ \ {\it Let $(M, g)$ be a complete non-compact Riemannian
$n$-manifold with $n \ge 2$ and $S_g \le 0$. Let $x_o \in M$ and let 
$B_R$ be the
open ball on $M$ with center at $x_o$ and radius equal to $R > 0$. If 
there exist
positive constants $C$ and $\alpha \in (0, n/2)$ such that}
$${\mbox {Vol}} (B_R) \le C R^\alpha$$ {\it for all $R > 0$, then $Q (M, 
g) \le
0$.}\\[0.2in]
{\bf Proof.} \ \ Assume that $Q (M, g) > 0$. Let $\alpha' \in (0, n/2)$ 
be a
positive constant such that $\alpha < \alpha'$. There exist positive 
constants
$c$  and $R_o$ such that
$$c^{2\over n} \le {{Q (M, g)}\over 4} \ \ \ \ {\mbox {and \ \ \ \ Vol}}
(B_R)
\le c R^{\alpha'}$$ for all $R \ge R_o$. For $R > 0$, let
$$\lambda_R = \inf \,\, \{ {{ \int_{B_R} |\bigtriangledown u |^2 dv_g 
}\over {
(\int_{B_R} | u |^2 dv_g ) }} \ | \ \ u \in C^\infty_o
(B_R)\,, u \not\equiv 0 \}\,.$$
Then
\begin{eqnarray*}
\lambda_R & \ge & \inf \,\, \{ {{ \int_{B_R} |\bigtriangledown u |^2 dv_g 
+ {{n -
2}\over {4 (n - 1)}} \int_{B_R} S_g u^2 dv_g }\over { (\int_{B_R} | u |^2 
dv_g )
}}
\ |
\ \ u \in C^\infty_o (B_R)\,, u \not\equiv 0 \}\\
& = & \inf \,\, \{ {{ (\int_{B_R} | u |^{{2n}\over {n - 2}} dv_g )^{{n - 
2}\over n}
}\over {\int_{B_R} | u |^2 dv_g}} \times {{ \int_M |\bigtriangledown u 
|^2 dv_g +
{{n - 2}\over {4 (n - 1)}} \int_M S_g u^2 dv_g }\over { (\int_M | u 
|^{{2n}\over
{n - 2}} dv_g )^{{n - 2}\over n} }}  \ | \ \ u \in C^\infty_o (B_R)\,, u
\not\equiv 0
\}\,.
\end{eqnarray*}
By H\"older's inequality we have
$$\int_{B_R} | u |^2 dv_g \le ( \int_{B_R}
|u|^{{2n}\over {n - 2}} dv_g )^{{n - 2}\over n} ( \int_{B_R} 1 \cdot dv_g
)^{2\over n}\,.$$ 
Therefore
$${{ (\int_{B_R} | u |^{{2n}\over {n - 2}} dv_g )^{{n - 2}\over n} }\over
{\int_{B_R} | u |^2 dv_g}} \ge {1\over {[ \mbox {Vol} (B_R) ]^{2\over 
n}}} \ge
{1\over {c^{2\over n} R^{{2\alpha'}\over n} }}\,, \ \ \ \ R \ge R_o\,.$$
As $Q (M, g) > 0$, we have 
$${{ \int_M |\bigtriangledown u |^2 dv_g + {{n - 2}\over {4 (n - 1)}} 
\int_M S_g
u^2 dv_g }\over { (\int_M | u |^{{2n}\over {n - 2}} dv_g )^{{n - 2}\over 
n} }}
> 0$$
for all $u \in C^\infty_o (M)\,, u \not\equiv 0$. 
Thus 
\begin{eqnarray*}
\lambda_R & \ge & {1\over {c^{2\over n} R^{{2\alpha'}\over n} }}   
\inf \,\, \{ {{ \int_M |\bigtriangledown u |^2 dv_g + {{n -
2}\over {4 (n - 1)}} \int_M S_g u^2 dv_g }\over { (\int_M | u 
|^{{2n}\over {n -
2}} dv_g )^{{n - 2}\over n} }} \ | \ \ u \in C^\infty_o (M)\,, u 
\not\equiv 0 \}
\\
& \ge & {4\over {R^{{2\alpha'}\over n} }}
\end{eqnarray*}
for $R \ge R_o$, as $c^{2\over n} \le Q (M, g)/4$. Let $k > R_o$ be an 
integer
and let
$\varphi_k \in C^\infty_o (B_{k + 1})$ be such that $\varphi_k (x) \ge 0$
for all $x \in B_{k + 1}$, $\varphi_k (x) = 1$ for
$x \in B_k$ and $|\bigtriangledown \varphi_k | \le 2$ [cf. 11]. We have
\begin{eqnarray*}
{\mbox {Vol\,\,}} (B_{k + 1}) - {\mbox {Vol\,\,}} (B_k) & \ge & {1\over 
4} \int_{B_{k +
1}} |\bigtriangledown \varphi_k |^2 dv_g\\
 & \ge & {1\over 4} \lambda_{k + 1} \int_{B_{k + 1}} | \varphi_k
|^2 dv_g\\
 & \ge & {{ {\mbox {Vol\,\,}} (B_k)}\over {(k + 1)^{{2\alpha'}\over n} }}\,.
\end{eqnarray*} 
Therefore
$${\mbox {Vol\,\,}} (B_{k + 1}) \ge (1 + {1\over {(k + 
1)^{{2\alpha'}\over n} }}
) {\mbox {Vol\,\,}} (B_k)$$
for all $k > R_o$. Let $\beta = 2 \alpha '/n < 1$. Given an integer $m >
n/2$, there exists an integer
$k_o > R_o$ such that for any integer $k > k_o$, we have
\begin{eqnarray*}
( {{k + 2}\over {k + 1}} )^m & = & (1 + {1\over {k + 1}})^m\\
& = & 1 + {{c_1
(m)}\over {k + 1}}  + {{c_2 (m)}\over {(k + 1)^2}} + \cdot \cdot \cdot + 
{{c_{m
- 1} (m)}\over {(k + 1)^{m - 1}}} + {1\over {(k + 1)^m}}\\ & \le & 1 + {{C
(m)}\over {k + 1}}\\ & \le & 1 + {1\over {(k + 1)^\beta}} \,,
\end{eqnarray*}
where $c_1 (m),..., c_{m - 1} (m)$ are positive constants and $C (m)$ is 
a large
positive constant. 
We
have
\begin{eqnarray*}
{\mbox {Vol\,\,}} (B_{k + 1}) & \ge & (1 + {1\over {(k + 1)^\beta}}) (1 +
{1\over k^\beta})
\cdot \cdot \cdot (1 + {1 \over {(k_o + 1)^\beta}}) {\mbox {Vol\,\,}} 
(B_{k_o})\\
& \ge & ({{k + 2}\over {k + 1}})^m \times ({{k + 1}\over {k }})^m \times 
\cdot
\cdot
\cdot
\times ({{k_o + 3}\over {k_o + 2}})^m \times ({{k_o +
2}\over {k_o + 1}})^m {\mbox {Vol\,\,}} (B_{k_o})\\
& = & (k + 2)^m {{ {\mbox {Vol\,\,}} (B_{k_o})}\over {(k_o + 1)^m}}
\end{eqnarray*}
for all $k > k_o$. This contradicts that ${\mbox {Vol\,\,}} (B_{k + 1}) 
\le c (k
+ 1)^{\alpha '}$ with $\alpha ' < n/2$, for all $k > k_o$. \ \ \ \ \ {\bf
Q.E.D.}\\[0.2in]
{\bf Corollary 2.2.} \ \ {\it Let $(M, g)$ be a complete Riemannian 
manifold of
non-positive scalar curvature. If the volume of $(M, g)$ is finite, then $Q
(M, g) \le 0$.}\\[0.2in]
{\bf Corollary 2.3.} \ \ {\it For $n \ge 3$, let $(M, g)$ be a simply 
connected,
complete, non-compact, conformally flat Riemannian manifold of non-positive
scalar curvature. Then for some positive constant $C$ and for any $\alpha 
< n/(n
+ 2)$,
${\mbox {Vol\,\,}} B_R
\ge C R^\alpha$ for all
$R$ large.}\\[0.2in]
{\bf Proof.} \ \ As    
$(M, g)$ is a simply connected, complete, non-compact, conformally flat
Riemannian manifold, we have [13]
$$Q (M, g) = Q (S_n, g_o) = {{n (n - 2) \omega_n^{2\over n}}\over 4}\,,$$
where $\omega_n$ denotes the volume of $(S^n, g_o)$. Let $C$ be a positive
constant with $C^{2\over n} \le Q (M, g)/4$. If the statement that Vol
$B_R
\ge C R^\alpha$ for all $R$ large does not hold for a positive constant 
$\alpha
< n/(n + 2)$, then there exists an increasing  sequence of positive numbers
$\{r_o, r_1,...., r_k,....\}$ such that 
$$\lim_{k \to \infty} r_k = \infty \ \ \ \ {\mbox {and}} \ \ \ \ \ {\mbox 
{Vol\,\,}}
B_{r_k}
\le C r_k^\alpha$$ 
for $k = 1, 2,....$ It follows from the proof of theorem 2.1 that for $k$
large,  
\begin{eqnarray*}
{\mbox {Vol\,\,}} B_{r_k} & \ge & ( 1 + {1\over { r_k^{{2 \alpha}\over n} 
}} )
{\mbox {Vol\,\,}} B_{r_k\, - 1}\\
& \cdot & \\
& \cdot & \\
& \cdot & \\
& \ge & ( 1 + {1\over { r_k^{{2 \alpha}\over n} }} )^{[r_k] - 1} {\mbox 
{Vol\,\,}}
B_1\\
& \ge & ( 1 + {{[r_k] - 1}\over { r_k^{{2 \alpha}\over n} }} ) {\mbox 
{Vol\,\,}}
B_1\\
& \ge & {{r_k}\over { r_k^{{2 \alpha}\over n} }}{\mbox {Vol\,\,}} B_1\\ 
& \ge & r_k^{\alpha '}{\mbox {Vol\,\,}} B_1\,,
\end{eqnarray*}
where $[r_k]$ is the integer part of $r_k$ and $\alpha '$ is a constant
such that $\alpha <
\alpha ' < n/(n + 2)$. When $k$ is large
enough, the last inequality contradicts that ${\mbox {Vol\,\,}} B_{r_k}
\le C r_k^\alpha$ for all non-negative integer $k$. \ \ \ \ \ {\bf
Q.E.D.}\\[0.2in]
{\bf Theorem 2.4.} \ \ {\it Let $(M, g)$ be a complete non-compact Riemannian
$n$-manifold with $n \ge 3$. Assume that $S_g \le 0$ and there exist 
constants
$C > 0$, $R_1 > 0$ and
$\beta \in (0, 1)$ such that 
$$S_g (x) \le - C d (x_o, x)^{-\beta}$$
for all $x \in M \setminus B_{R_1}$, where $x_o
\in M$ is a fixed point and $d (x_o, x)$ is the distance between $x_o$
and
$x$ in
$(M, g)$. If there exist a positive constant $C'$ and a positive integer $m$
such that}
$${\mbox {Vol}} (B_R) \le C' R^m$$ 
{\it for all $R > 0$, then $Q (M, g)
< 0$.}\\[0.2in] 
{\bf Proof.} \ \ Given a positive number $\delta$, let $u_o \in 
C^\infty_o (M)$
be a smooth function such that 
$${{4 (n - 1) }\over {n - 2}} \Delta u_o = \delta \ \ \ {\mbox {on}} \ \ \
B_{R_1}\,.$$
Assume that the support of $u_o$ is inside $B_{R_2}$ for some constant 
$R_2 >
R_1$. We can find a positive constant $c_o$ such that 
$$S_g (x) (c_o + u_o (x)) - {{4 (n - 1) }\over {n - 2}} \Delta u_o (x) \le
- \delta$$ 
for all $x \in B_{R_2}$. This is because $S_g \le 0$ and  
$$S_g (x) \le - C R_2^{-\beta}$$
for all $x \in B_{R_2} \setminus B_{R_1}$. Let $u = c_o + u_o$. Then $u$ 
is a
smooth positive function on $M$. Furthermore, $u (x) = c_o$ for all $ x 
\in M
\setminus B_{R_2}$. Let $g' = u^{4/(n - 2)} g$ and let $S'$ be the scalar
curvature of the metric $g'$. Then [9]
$$S' (x) = u^{ - {{n + 2}\over {n - 2}} } [S_g (x) u (x) - {{4 (n - 1) 
}\over {n
- 2}} \Delta u_o (x) ] \le - \epsilon$$
for all $x \in B_{R_2}$, where $-\epsilon$ is a suitable constant. And
$$S' (x) = c_o^{ - {{4}\over {n - 2}} } S_g (x)$$ 
for all $x \in  M \setminus B_{R_2}$. There exist positive constants 
$c_1$ and
$C_1$ such that $c_1 \le u \le C_1$ on $M$. Therefore
$$c_1^{- {2 \over {n - 2}} } d' (x_o, x) \le d (x_o, x)  \le  C_1^{2
\over {n - 2}} d' (x_o, x)$$ for all $x \in M$, where $d' (x_o, x)$ is the
distance between
$x_o$ and
$x$ with respect to the Riemannian metric $g'$. Thus we can find a positive
constant $R'_1$ such that 
$$S' (x) \le  - C'' d' (x_o, x)^{-\beta}$$
for all $x \in M$ with $d' (x_o, x) > R'_1$, where $C''$ is a positive
constant. And for $d' (x_o, x) \le R'_1$, we can find a positive constant
$\epsilon'$ such that
$$S' (x) \le - \epsilon'\,.$$
Hence we can find a positive number $R_o$ such that for all $R > R_o$, we 
have
$$S' (x) \le  - C_o R^{-\beta}$$
for all $x \in M$ with $d' (x_o, x) \le R$, where $C_o$ is a positive 
integer. 
As $Q (M, g') = Q (M, g)$, by above, we may assume without loss of
generality that for all $R > R_o$,  
$$S_g (x) \le  - C_o R^{-\beta}$$ 
for all $x \in M$ with $d (x_o, x) \le  R\,.$\\[0.02in]
\hspace*{0.5in}Assume that $Q (M, g) \ge 0$. Take $\beta' \in (0, 1)$ such
that $\beta' > \beta$. Then we can find a positive number $R_3 > R_o$ such
that for all $R > R_3$,  
$$S_g (x) \le - {{4^2 (n - 1)}\over {n - 2}} \, {1\over {R^{\beta'} }}$$
for all $x \in M$ with $d (x_o, x) \le R$. For $R > R_3$, we have
\begin{eqnarray*}
& \ & \lambda_R \\
& \ge & \inf \,\, \{ {{ \int_{B_R} |\bigtriangledown u |^2 dv_g +
{{n - 2}\over {4 (n - 1)}} \int_{B_R} S_g u^2 dv_g }\over { (\int_{B_R} | 
u |^2
dv_g ) }} - {{ {{n - 2}\over {4 (n - 1)}} \int_{B_R} S_g u^2 dv_g }\over 
{(\int_{B_R} | u |^2 dv_g ) }}
\ |
\ \ u \in C^\infty_o (B_R)\,, u \not\equiv 0 \}\\ 
& \ge & \inf \,\, \{ {{ (\int_{B_R} | u |^{{2n}\over {n - 2}} dv_g )^{{n -
2}\over n} }\over {\int_{B_R} | u |^2 dv_g}} \times {{\int_{B_R}
|\bigtriangledown u |^2 dv_g + {{n - 2}\over {4 (n - 1)}} \int_{B_R} S_g u^2
dv_g }\over { (\int_{B_R} | u |^{{2n}\over {n - 2}} dv_g )^{{n - 2}\over 
n} }}
\ |
\ \ u \in C^\infty_o (B_R)\,, u \not\equiv 0 \}\\
& \ &  \ \ \ \ \ + {4 \over {R^{\beta'} }}\\
& \ge & {4 \over {R^{\beta'} }}\,,
\end{eqnarray*}
as $Q (M, g) \ge 0$ implies that 
$${{\int_{B_R} |\bigtriangledown u |^2 dv_g + {{n - 2}\over {4 (n - 1)}}
\int_{B_R} S_g u^2 dv_g }\over { (\int_{B_R} | u |^{{2n}\over {n - 2}} dv_g
)^{{n - 2}\over n} }} \ge 0$$
for all $u \in C^\infty_o (B_R)\,, u \not\equiv 0$. It follows as in 
theorem 2.1
that for any positive integer
$m' > m$, there exists an integer $k_o > R_3$ such that for all $k > 
k_o$, we
have
$${\mbox {Vol\,\,}} (B_{k + 1}) \ge  (k + 2)^{m'} {{ {\mbox {Vol\,\,}}
(B_{k_o})}\over {(k_o + 1)^{m'}}}$$
for all $k > k_o$. This contradicts that ${\mbox {Vol\,\,}}
(B_{k + 1}) \le C' R^m$ with $m < m'$, for all $k > k_o$.
\ \ \ \ \ {\bf Q.E.D.}\\[0.2in]
{\bf Corollary 2.5.} \ \ {\it Let $(M, g)$ be a simple connected, 
complete, 
non-compact conformally flat Riemannian
$n$-manifold with $n \ge 3$. Assume that $S_g \le 0$ and there exist 
constants
$C > 0$,
$R_1 > 0$ and
$\beta \in (0, 1)$ such that 
$$S_g (x) \le - C d (x_o, x)^{-\beta}$$ for all $x \in M \setminus B_{R_1}$,
where
$x_o
\in M$ is a fixed point. Then for any positive numbers $m$ and $C'$, 
there exists
a positive constant
$R_1$ such that}
$${\mbox {Vol\,\,}} (B_R) \ge C' R^m$$  
{\it for all $R > R_1$. }\\[0.2in] 
{\bf Remark.} \ \ Corollary 2.5 shows that on the  Euclidean space ${\bf
R}^n$, $n \ge 3$, there does not exist a Riemannian metric $g'$ which is 
uniformally equivalent and conformal to the Euclidean metric and with
non-positive scalar curvature $S'$ which satisfies
$$S' (x) \le - C |x|^{- \beta}$$
for some constants $\beta \in (0, 1)$ and $C > 0$, for all $x$ outside a
compact domain containing the origin. Such a metric exists if the scalar
curvature is allowed to drop to zero faster. In fact, Ni [11] has 
constructed a
metric
$g'$ in
${\bf R}^n$ which is uniformally equivalent and conformal to the Euclidean
metric, with $S' < 0$ and for all $x$ outside a compact domain containing the
origin, 
$$S' (x) \le -C |x|^{- l}$$
for some constant $l > 2$ and $C > 0$.\\[0.2in] 
{\bf Theorem 2.6.} \ \ {\it {Let $(M, g)$ be a compact conformally flat
$n$-manifold with $n \ge 3$ and scalar curvature $S_g \leq 0$. Then the
fundamental group of
$M$ has exponential growth unless $S_g \equiv 0.$}}\\[0.2in]
{\bf Proof.} \ \ If $S_g \not\equiv 0$, then $S_g$ is negative somewhere and
hence the conformal invariant $Q (M, g)$ is negative. By a conformal 
change of
the metric
$g$, we may assume that $S_g \equiv - c^2$ with $c > 0$ being a constant 
[9]. Let
$\tilde M$ be the universal covering of $M$, equipped with the pull back 
metric.
Then
$$Q (\tilde M, g) = Q (S_n, g_o) = {{n (n - 2) \omega_n^{2\over n}}\over 
4}\,.$$
Given a point $p
\in \tilde M$ and $R > 0$, let $B_R$ be the ball in $\tilde M$ with center
at $p$ and radius $R$. As in the proof of theorem 2.4, we have 
$$\lambda_R \ge {{n - 2}\over {4 (n - 1)}} c^2$$ 
and hence there exist positive constants $C$ and
$\delta$ such that (c.f. proposition 1.2 of [13]) 
$${\rm Vol\,\,} (B_R) \ge C e^{\delta R}$$
for all $R$ large. Then an argument as in [10] shows
that $\pi_1 (M)$ has exponential growth. \ \ \ \ \ {\bf Q.E.D.}\\[0.2in]
{\bf Theorem 2.7.} \ \ {\it For $n \ge 3$, let $(M, g)$ be a compact 
conformally
flat
$n$-manifold  of non-positive scalar curvature. If $\pi_1 (M)$ is almost
solvable, then $(M, g)$ is a flat manifold.}\\[0.2in]
{\bf Proof.} \ \ Since $\pi_1 (M)$ is almost solvable, the holonomy group of
$(M, g)$ is also almost solvable. By a result in [6],
$(M, g)$ is covered by a conformally flat manifold which is either 
conformally
diffeomorphic to the
$n$-sphere
$S^n$, a flat
$n$-torus, or a Hopf manifold
$S^1 \times S^{n-1}$. If $M$
is covered by $S^1 \times S^{n-1}$, then  $S^1 \times S^{n-1}$ admits a
conformally flat metric of non-positive scalar curvature. Since 
$\pi_1 (S^1 \times S^{n-1})$ does not have exponential growth, theorem
2.6 implies that the scalar curvature of $\bar g$ is identically equal to 
zero.
Then by  theorem 4.5 in [13], $(S^1 \times S^{n-1}\,, \bar g)$ is conformally
equivalent to 
$\Omega / \Gamma$, where $\Omega$ is the image of a developing map and 
$\Gamma$
is the corresponding holonomy group. Furthermore, $\Omega$ is the domain of
discontinuity of $\Gamma$.  Since $\Gamma \cong {\bf Z}$, by using a 
result in
[6], we have $\Gamma$ is conjugate to a cyclic  group of similarities. Hence
$\Omega / \Gamma$ is conformally equivalent to $S^1 \times S^{n-1}$ with 
the 
standard product metric, denoted by $h_1$. If we denote 
$$\psi : (S^1 \times S^{n-1}\,, \bar g) \to (S^1 \times S^{n-1}\,, h_1)$$ 
a conformal equivalence and $h^* = \psi^* h_1$, then there exists a 
positive 
function $u$ such that $\bar g = u^{{4}\over {n-2}} h^*$. Hence
$$\Delta_{\bar g}\, u - c_n S_{\bar g}\, u = -c_n S_{h^{*}} u^{{{n+2}\over
{n-2}}}\,,$$
where $S_{\bar g}$ and $S_{h^*}$ are the scalar curvature of $\bar
g$ and
$h^{*}$, respectively. Now $S_{\bar g} = 0$ and $S_{h^*} > 0$, the maximum
principle implies that 
$u$ is a constant, which is impossible.\\[0.02in]
\hspace*{0.5in}Suppose that $M$ is covered by $S^n$. By pulling back the 
metric
$g$ to $S^n$, we obtain a metric $g'$. Then $g'$ is a conformally flat 
metric on
$S^n$ with non-positive scalar curvature.  Let $g_o$ be the standard 
metric on
$S^n$.  The developing map
$$\phi: (S^n, g') \rightarrow (S^n, g_o)$$  is a conformal equivalence. Let
$g^{*} = \phi^{*} g_o$.  Hence there exists a function $u > 0$ such that 
$g^{*}
= u^2 g'$ and 
$$\Delta_{g'}\, u - c_n S_{g'}\, u = -c_n S_{g^{*}} u^{{{n+2}\over 
{n-2}}}\,,$$\\
where $S_{g'}$ and $S_{g^{*}}$ are the scalar curvature of $g'$ and
$g^{*}$, respectively. In fact, $S_{g^{*}} = n (n-1)$ and
$S_{g'} \le 0$.  The curvature conditions imply ${\Delta}_{g'}\, u \le 0$.
Therefore $u$ is a constant, which is impossible. So we conclude that
$M$ is covered by a flat torus $T^n$. This implies that $\pi_1 (M)$ is almost
nilpotent and has polynomial growth [3]. Theorem 2.6 shows that the scalar
curvature of $g$ is identically equal to zero. By pulling back the metric 
$g$ to
$T^n$, we have a metric on $T^n$ which has zero scalar curvature. But any 
scalar
flat metric on $T^n$ is flat [8]. So the pull back metric is flat, hence $(M,
g)$ is flat. \ \ \ \ \ {\bf Q.E.D.}\\[0.2in]
{\bf Theorem 2.8.} \ \ {\it For $n = 4, 6$, let $M$ be a compact, 
orientable  
$n$-manifold with zero Euler characteristic. If $\pi_1 (M)$ does not have
exponential growth, then any conformally flat metric on $M$ with non-positive
scalar curvature is flat.}\\[0.2in]
{\it Proof.} \ \  Let $g$ be a conformally flat metric on $M$ with 
non-positive 
scalar curvature. For $n = 4$, let $B$ be the Gauss-Bonnet integrand of $(M,
g)$   and let $\cal R$ and $Ric$ be the curvature tensor and Ricci tensor 
of $g$,
respectively. Then there exist universal  constants [3] $\alpha,\  \beta 
\ {\rm
and} \
\gamma$ such that                                                        
$$B=\alpha {\mid \cal R \mid}^2 + \beta{\mid Ric \mid}^2 + \gamma{\mid S_g
\mid}^2\,.$$\\ Here $R$ is the scalar curvature of $g$. Since $(M, g)$ is
conformally flat, we can write  ${\mid
\cal R \mid}^2$ in terms of the other two norms. That is, 
$B = a{\mid Ric \mid}^2 + b{\mid S_g \mid}^2$. Evaluate $a$ and $b$ on $S^1
\times S^3$ with the standard product metric and $S^4$, we 
obtain              
$$a=-{{1}\over {6{\rm Vol\,\,}(S^4)}}\,,$$
and                                               
$$b= {1\over {18 {\rm Vol\,\,} (S^{4})}}\,.$$  Using the Gauss-Bonnet 
theorem, we
obtain                                    
$$0 = \chi (M)= {{-1}\over {6 {\rm Vol\,\,} (S^{4})}} \int_M {\mid Ric 
\mid}^2
dv_g +  {1\over {18 {\rm Vol\,\,} (S^{4})}}\int_M |S_g|^2 dv_g\,.$$\\ Since
$\pi_1 (M)$ does not have exponential growth, theorem 2.6 implies that $S_g
\equiv 0$. The above formula gives
$Ric \equiv 0$. Since
$g$ is conformally flat, the Weyl tensor is identically equal to zero. 
Therefore
$g$ is flat.\\[0.02in]
\hspace*{0.5in}For $n = 6$, we use equation (2.4) in [12] to obtain 
$$384\pi^3 \chi (M) = - \int_M |\bigtriangledown {\cal R}|^2 dv_g\,.$$\\
Since $\chi (M) = 0$, therefore $(M, g)$ is locally  symmetric. Then the
classification theorem of conformally flat symmetric spaces in  dimension six
[7] implies that the universal covering of 
$(M, g)$ is isometric to one of the following symmetric spaces:
$${\bf R}^6\,, \ \ S^6 (c)\,, \ \ {\bf H}^6 (-c)\,, \ \ {\bf R} \times S^5
(c)\,, \ \ {\bf R} \times {\bf H}^5 (-c)\,, \ \ S^2 (c) \times {\bf H}^4
(-c)\,,$$
$$S^4 (c) \times {\bf H}^2 (c)\,, \ \  S^3 (c) \times {\bf H}^3 (-c)\,,$$\\
where ${\bf H}^n (-c)$ is the $n$-dimensional simply connected complete 
manifold
of constant  sectional curvature equal to $-c$ and $S^n (c)$ is the 
sphere with
sectional  curvature equal to $c$. Since $(M, g)$ has non-positive scalar
curvature, it cannot be covered by $S^6 (c)$, ${\bf R} \times S^5 (c)$ or 
$S^4 (c) \times {\bf H}^2 (c)$. If $M$ were covered by ${\bf H}^6 (-c)$, 
${\bf R} \times {\bf H}^5 (-c)$, or $S^3 (c) \times {\bf H}^3 (-c)$, then 
$\pi_1
(M)$ would have  exponential growth, but this is not true. Hence $(M, g)$ can
only be covered by ${\bf R}^6$ and so $g$ is flat. \ \ \ \ \ {\bf Q.E.D.}\\

\pagebreak

\vspace{0.5in}

\centerline{\bf \Large References}

\vspace{0.3in}

[1] Avez, A.: {\it {Varietes Riemanniennes sans points focaux}}, C.R. 
Acad. Sc.
Paris {\bf {270}} (1970), 188-191.\\[0.1in]
[2] Ballmann, W., Gromov, M., Schroeder, V.: {\it Manifolds of nonpositive
curvature}, Birkh\"auser, 1985.\\[0.1in]
[3] Berger, M., Gauduchon, P., Mazet, E.: {\it {Le Spectre d'une Variete
Riemanniene}}, Lecture Notes in Mathematics, Vol. 194, Springer-Verlag, New
York.\\[0.1in]
[4] Cheeger, H., Ebin, D.: {\it Comparison theorems in Riemannian geometry},
North-Holland, 1975.\\[0.1in]
[5] Gromov, M.: {\it {Groups of polynomial growth and expanding maps,}} IHES
Publications Math\'ematiques {\bf {53}} (1981), 53-73.\\[0.1in]
[6] Kamishima, Y.: {\it Conformally flat manifolds whose development maps 
are not
surjective {\rm I}}, Transactions of the AMS {\bf 294} (1986), 
607-623.\\[0.1in]
[7] Kurita, M.: {\it On the holonomy group of the conformally flat Riemannian
manifold,}  Nagoya Math. J. {\bf 9} (1955), 161-171.\\[0.1in]
[8] Lawson, H. Jr., Michelsohn, M.: {\it {Spin Geometry,}} Princeton 
Mathematical
Series: 38, Princeton University Press, 1989.\\[0.1in]
[9] Lee, J., Parker, T.: {\it The Yamabe Problem}, Bulletin of the AMS  
{\bf 17}
(1987), 37-91.\\[0.1in]
[10] Milnor, J.: {\it {A note on curvature and fundamental group}}, J. 
Diff. Geom.
{\bf {2}} (1968), 1-7.\\[0.1in]
[11] Ni, W. M.: {\it On the elliptic equation $\Delta u + K (x) u^{{n + 
2}\over
{n - 2}} = 0$}, Indiana Univ. Math. J. {\bf 31} (1982), 493-529.\\[0.1in]
[12] Perrone, D.: {\it Osservazioni sulla caratteristica di
Eulero-Poincar\'e' di varieta'  Riemanniane conformemente piatte,} Note di
Matematica {\bf 3} (1983), 173-181.\\ [0.1in]
[13] Schoen, R., Yau, S.-T.: {\it {Conformally flat manifolds, Kleinian 
groups
and scalar curvature}}, Invent. Math. {\bf {92}} (1988), 47-71.

\end{document}